\documentclass[12pt]{article}
\usepackage{amssymb,amsmath,amscd, amsthm,stmaryrd,verbatim, mathrsfs}

\begin{document}
\baselineskip=15pt

\newcommand{\la}{\langle}
\newcommand{\ra}{\rangle}
\newcommand{\psp}{\vspace{0.4cm}}
\newcommand{\pse}{\vspace{0.2cm}}
\newcommand{\ptl}{\partial}
\newcommand{\dlt}{\delta}
\newcommand{\sgm}{\sigma}
\newcommand{\al}{\alpha}
\newcommand{\be}{\beta}
\newcommand{\G}{\Gamma}
\newcommand{\gm}{\gamma}
\newcommand{\vs}{\varsigma}
\newcommand{\Lmd}{\Lambda}
\newcommand{\lmd}{\lambda}
\newcommand{\td}{\tilde}
\newcommand{\vf}{\varphi}
\newcommand{\yt}{Y^{\nu}}
\newcommand{\wt}{\mbox{wt}\:}
\newcommand{\rd}{\mbox{Res}}
\newcommand{\ad}{\mbox{ad}}
\newcommand{\stl}{\stackrel}
\newcommand{\ol}{\overline}
\newcommand{\ul}{\underline}
\newcommand{\es}{\epsilon}
\newcommand{\dmd}{\diamond}
\newcommand{\clt}{\clubsuit}
\newcommand{\vt}{\vartheta}
\newcommand{\ves}{\varepsilon}
\newcommand{\dg}{\dagger}
\newcommand{\tr}{\mbox{Tr}}
\newcommand{\ga}{{\cal G}({\cal A})}
\newcommand{\hga}{\hat{\cal G}({\cal A})}
\newcommand{\Edo}{\mbox{End}\:}
\newcommand{\for}{\mbox{for}}
\newcommand{\kn}{\mbox{ker}}
\newcommand{\Dlt}{\Delta}
\newcommand{\rad}{\mbox{Rad}}
\newcommand{\rta}{\rightarrow}
\newcommand{\mbb}{\mathbb}
\newcommand{\lra}{\Longrightarrow}
\newcommand{\X}{{\cal X}}
\newcommand{\Y}{{\cal Y}}
\newcommand{\Z}{{\cal Z}}
\newcommand{\U}{{\cal U}}
\newcommand{\V}{{\cal V}}
\newcommand{\W}{{\cal W}}
\newcommand{\sta}{\theta}
\setlength{\unitlength}{3pt}
\newcommand{\msr}{\mathscr}

\begin{center}{\Large \bf Conformal Oscillator Representations \\ of Orthogonal Lie Algebras} 
\end{center}
\vspace{0.2cm}

\begin{center}{\large Xiaoping Xu
\footnote{Research supported
 by NSFC Grants 11171324 and  11321101.}}\end{center}
\begin{center}{Hua Loo-Keng Key Mathematical Laboratory\\
Institute of Mathematics, Academy of Mathematics \& System
Sciences\\ Chinese Academy of Sciences, Beijing 100190, P.R. China
}\end{center}

\begin {abstract}
\quad

The conformal transformations with respect to the metric defining
the orthogonal Lie algebra $o(n,\mbb C)$ give rise to a
one-parameter ($c$) family  of inhomogeneous first-order
differential operator representations of the orthogonal Lie algebra
$o(n+2,\mbb C)$. Letting these operators act on the space of
exponential-polynomial functions that depend on  a parametric vector
$\vec a\in \mbb C^n$, we prove that the space forms an irreducible
$o(n+2,\mbb C)$-module for any $c\in\mbb C$ if $\vec a$ is not on a
certain hypersurface. By partially swapping differential operators
and multiplication operators, we obtain more general differential
operator representations of $o(n+2,\mbb C)$ on the polynomial
algebra $\msr C$ in $n$ variables. Moreover, we prove that $\msr C$
forms an infinite-dimensional irreducible weight $o(n+2,\mbb
C)$-module with finite-dimensional weight subspaces if $c\not\in\mbb
Z/2$. \vspace{0.3cm}

\noindent{\it Keywords}:\hspace{0.3cm} orthogonal Lie algebra;
differential operator; oscillator
 representation; irreducible module; polynomial algebra; exponential-polynomial function.

\end{abstract}

\section {Introduction}

\quad$\;$ A module of a finite-dimensional simple Lie algebra is
called a {\it weight module} if it is a direct sum of its weight
subspaces. A module of a finite-dimensional simple Lie algebra is
called {\it cuspidal} if it is not induced from its proper parabolic
subalgebras. Infinite-dimensional irreducible weight modules of
finite-dimensional simple Lie algebras with finite-dimensional
weight subspaces have been intensively studied by the authors in
[BBL], [BFL], [BHL], [BL1], [BL2], [Fs], [Fv], [M]. In particular,
Fernando [Fs] proved that such modules must be cuspidal or
parabolically induced. Moreover, such cuspidal modules exist only
for special linear Lie algebras and symplectic Lie algebras. A
similar result was independently obtained by Futorny [Fv]. Mathieu
[M] proved that these cuspidal such modules
 are irreducible components in the tensor
modules of their multiplicity-free modules with finite-dimensional
modules. Although the structures of irreducible weight modules of
finite-dimensional simple Lie algebra with finite-dimensional weight
subspaces were essentially determined by Fernando's result in [Fs]
and Methieu's result in [M], explicit structures of such modules are
not that known. It is important to find explicit natural
realizations of them.

The $n$-dimensional conformal group with respect to Euclidean metric
$(\cdot,\cdot)$ is generated by the translations, rotations,
dilations and special conformal transformations
$$\vec x\mapsto\frac{\vec x-(\vec x,\vec x)\vec b}{(\vec b,\vec b)
(\vec x,\vec x)-2(\vec b,\vec x)+1}.\eqno(1.1)$$ Conformal groups
play important roles in geometry, partial differential equations and
quantum physics. The conformal transformations with respect to the
metric defining $o(n,\mbb{C})$ give rise to an inhomogeneous
representation of the Lie algebra $o(n+2,\mbb{C})$ on the polynomial
algebra in $n$ variables. Using Shen's mixed product for Witt
algebras in [S] and the above representation, Zhao and the author
[XZ]
 constructed   a new functor from $o(n,\mbb{C})$-{\bf Mod} to $o(n+2,\mbb{C})$-{\bf Mod} and derived a condition
 the functor to map a finite-dimensional irreducible
$o(n,\mbb{C})$-module to an infinite-dimensional irreducible
$o(n+2,\mbb{C})$-module. Our general frame also gave a direct
polynomial extension from irreducible $o(n,\mbb{C})$-modules to
irreducible $o(n+2,\mbb{C})$-modules.

The work [XZ] lead to a one-parameter ($c$) family of inhomogeneous
first-order differential operator (oscillator) representations of
$o(n+2,\mbb{C})$. Letting these operators act on the space of
exponential-polynomial functions that depend on  a parametric vector
$\vec a\in \mbb C^n$, we prove in this paper that the space forms an
irreducible $o(n+2,\mbb C)$-module for any $c\in\mbb C$ if $\vec a$
is not on a certain hypersurface. By partially swapping differential
operators and multiplication operators, we obtain more general
differential operator (oscillator) representations of $o(n+2,\mbb
C)$ on the polynomial algebra $\msr C$ in $n$ variables. Moreover,
we prove that $\msr C$ forms an infinite-dimensional irreducible
weight $o(n+2,\mbb C)$-module with finite-dimensional weight
subspaces if $c\not\in\mbb Z/2$.
 Our results are
extensions of Howe's oscillator construction of infinite-dimensional
multiplicity-free irreducible representations for $sl(n,\mbb{C})$
(cf. [H]).

For any two integers $p\leq q$, we denote
$\ol{p,q}=\{p,p+1,\cdots,q\}$. Let $E_{r,s}$ be the square matrix
with 1 as its $(r,s)$-entry and 0 as the others. Fix a positive
integer $n$.  Denote
$$A_{i,j}=E_{i,j}-E_{n+1+j,n+1+i},\;\;B_{i,j}=E_{i,n+1+j}-E_{j,n+1+i},\;\;C_{i,j}=E_{n+1+i,j}-E_{n+1+j,i}\eqno(1.2)$$
for $i,j\in\ol{1,n+1}$. Then the split even orthogonal Lie algebra
$$ o(2n+2,\mbb{C})=\sum_{i,j=1}^{n+1}
(\mbb{C}A_{i,j}+\mbb{C}B_{i,j}+\mbb{C}C_{i,j}).\eqno(1.3)$$  Set
$$D=\sum_{r=1}^nx_r\partial_{x_r}+\sum_{s=1}^ny_s\ptl_{y_s},\;\;\eta=\sum_{i=1}^nx_iy_i.\eqno(1.4)$$ According
to Zhao and the author's work [XZ],  we have the following
one-parameter generalization $\pi_c$ of the conformal representation
of $o(2n+2,\mbb{C})$:
$$\pi_c(A_{i,j})=x_i\ptl_{x_j}-y_j\ptl_{x_i},\;\pi_c(B_{i,j})=x_i\ptl_{y_j}-x_j\ptl_{y_i},\;
\pi_c(C_{i,j})=y_i\ptl_{x_j}-y_j\ptl_{x_i},\eqno(1.5)$$
$$\pi_c(A_{n+1,i})=\ptl_{x_i},\;\;\pi_c(A_{n+1,n+1})=-D-c,\;\;\pi_c(B_{i,n+1})=-\ptl_{y_i},\eqno(1.6)$$
$$\pi_c(A_{i,n+1})=\eta\ptl_{y_i}-x_i(D+c),\;\;\pi_c(C_{n+1,i})=y_i(D+c)-\eta\ptl_{x_i}\eqno(1.7)$$
for $i,j\in\ol{1,n}$. For $\vec a=(a_1,a_2,...,a_n)^t,\;\vec
b=(b_1,b_2,...,b_n)^t\in\mbb{C}^n$, we put
$$\vec a\cdot\vec
x=\sum_{i=1}^na_ix_i,\qquad\vec b\cdot\vec
y=\sum_{i=1}^nb_iy_i.\eqno(1.8)$$ Let ${\msr
A}=\mbb{C}[x_1,...,x_n,y_1,...,y_n]$ be the algebra of polynomials
in $x_1,...,x_n,y_1,...,y_n$. Moreover, we set
$${\msr A}_{\vec a,\vec b}=\{fe^{\vec a\cdot\vec
x+\vec b\cdot\vec y}\mid f\in{\msr A}\}.\eqno(1.9)$$
 Denote by $\pi_{c,\vec a,\vec b}$ the
representation $\pi_c$ of $o(2n+2,\mbb C)$ on $\msr A_{\vec a,\vec
b}$.

Fix $n_1,n_2\in\ol{1,n}$ with $n_1\leq n_2$.
 Changing operators $\ptl_{x_r}\mapsto -x_r,\;
 x_r\mapsto
\ptl_{x_r}$  for $r\in\ol{1,n_1}$ and $\ptl_{y_s}\mapsto -y_s,\;
 y_s\mapsto\ptl_{y_s}$  for $s\in\ol{n_2+1,n}$ in the
 representation  $\pi_c$ of $o(2n+2,\mbb{C})$, we get another
differential-operator representation $\pi_c^{n_1,n_2}$ of
$o(2n+2,\mbb{C})$ on $\msr A$. We call $\pi_c$ and $\pi_c^{n_1,n_2}$
the {\it conformal oscillator representations of $o(2n+2,\mbb{C})$}
in terms of physics terminology.
  In this paper, we prove:\psp

{\bf Theorem 1}. {\it The representation $\pi_{c,\vec a,\vec b}$ of
$o(2n+2,\mbb{C})$   is irreducible for any $c\in\mbb{C}$ if
$\sum_{i=1}^na_ib_i\neq 0$. Moreover, the representation
$\pi_c^{n_1,n_2}$ of $o(2n+2,\mbb{C})$ is irreducible for any
$c\in\mbb{C}\setminus(\mbb Z/2)$, and its underlying module ${\msr
A}$ is an infinite-dimensional irreducible weight
$o(2n+2,\mbb{C})$-module with finite-dimensional weight subspaces.
}\psp

Set
$$K_i=E_{0,i}-E_{n+i+1,0},\qquad
K_{n+1+i}=E_{0,n+1+i}-E_{i,0}\qquad\for\;\;i\in\ol{1,n+1}.\eqno(1.10)$$
Then the split odd orthogonal Lie algebra
$$o(2n+3,\mbb{C})= o(2n+2,\mbb{C})+\sum_{i=1}^{2n+2}\mbb{C}K_i.\eqno(1.11)$$
Moreover, we redefine
$$D=\sum_{r=0}^nx_r\partial_{x_r}+\sum_{r=1}^ny_r\ptl_{y_r}\qquad\eta=\frac{1}{2}x_0^2+\sum_{i=1}^nx_iy_i.
\eqno(1.12)$$
 According
to Zhao and the author's work [XZ],  we have the following
one-parameter generalization of the conformal representation $\pi_c$
of $o(2n+3,\mbb{C})$: $\pi_c|_{o(2n+2,\mbb{C})}$ is given in
(1.5)-(1.7) with $D$ and $\eta$ in (1.12),
$$\pi_c(K_i)=x_0\ptl_{x_i}-y_i\ptl_{x_0},\;\;\pi(K_{n+1+i})=x_0\ptl_{y_i}-x_i\ptl_{x_0}\qquad\for\;\;i\in\ol{1,n},
\eqno(1.13)$$
$$\pi_c(K_{n+1})=x_0(D+c)-\eta\ptl_{x_0},\qquad \pi_c(K_{2n+2})=-\ptl_{x_0}.\eqno(1.14)$$
Fix $n_1,n_2\in\ol{1,n}$ with $n_1\leq n_2$.
 Changing operators $\ptl_{x_r}\mapsto -x_r,\;
 x_r\mapsto
\ptl_{x_r}$  for $r\in\ol{1,n_1}$ and $\ptl_{y_s}\mapsto -y_s,\;
 y_s\mapsto\ptl_{y_s}$  for $s\in\ol{n_2+1,n}$ in the above
 representation of $o(2n+3,\mbb{C})$, we get another
differential-operator representation $\pi_c^{n_1,n_2}$ of
$o(2n+3,\mbb{C})$. Again call the representations $\pi_c$ and
$\pi_c^{n_1,n_2}$ of $o(2n+3,\mbb{C})$ {\it conformal oscillator
representations} in terms of physics terminology.

Let ${\msr B}=\mbb{C}[x_0,x_1,...,x_n,y_1,...,y_n]$ be the algebra
of polynomials in $x_0,x_1,...,x_n,y_1,...,y_n$. Redenote
$$\vec a\cdot\vec x=\sum_{i=0}^na_ix_i\qquad\for\;\;\vec
a=(a_0,a_1,...,a_n)^t\in\mbb{C}^{1+n}.\eqno(1.15)$$ Fix $\vec a
\in\mbb{C}^{1+n},\; \vec b\in\mbb C^n$ and $n_1,n_2\in\ol{1,n}$ with
$n_1\leq n_2$.  We set
$${\msr B}_{\vec a,\vec b}=\{fe^{\vec a\cdot\vec
x+\vec b\cdot\vec y}\mid f\in{\msr B}\}\eqno(1.16)$$ (cf. (1.8)).
Denote by $\pi_{c,\vec a,\vec b}$ the representation $\pi_c$ of
$o(2n+3,\mbb{C})$ on ${\msr B}_{\vec a,\vec b}$.

In [XZ], Zhao and the author proved that the representation
$\pi_{c,\vec 0,\vec 0}$ of $o(2n+3,\mbb{C})$ is irreducible if and
only if $c\not\in -\mbb{N}$. The following is our second main
theorem in this paper.\psp

{\bf Theorem 2}. {\it The  representation $\pi_{c,\vec a,\vec b}$
 of $o(2n+3,\mbb{C})$ is irreducible for any $c\in\mbb C$ if
 $a_0^2+2\sum_{i=1}^na_ib_i\neq 0$. Moreover, the representation $\pi_c^{n_1,n_2}$ of
$o(2n+3,\mbb{C})$ is irreducible for any $c\in\mbb{C}\setminus(\mbb
Z/2)$, and its underlying module ${\msr B}$ is an
infinite-dimensional irreducible weight $o(2n+3,\mbb{C})$-module
with finite-dimensional weight subspaces. }\psp

In Section 2, we prove Theorem 1. The proof of Theorem 2 is given in
Section 3.

\section{Proof of Theorem 1}

First we want to prove:\psp

 {\bf Theorem 2.1}. {\it The representation $\pi_{c,\vec a,\vec
b}$ of $o(2n+2,\mbb{C})$ is irreducible if $\sum_{i=1}^na_ib_i\neq
0$ for any $c\in\mbb{C}$.}

{\it Proof}. By symmetry, we may assume $a_1\neq 0$. Let ${\msr M}$
be a nonzero $o(2n+2,\mbb{C})$-submodule of ${\msr A}_{\vec a,\vec
b}$. Take any $0\neq fe^{\vec a\cdot\vec x+\vec b\cdot\vec y}\in
\msr{M}$ with $f\in \msr{A}$.  Let $\msr{A}_k$ be the subspace of
homogeneous polynomials with degree $k$. Set
$$\msr{A}_{\vec a,\vec b,k}=\msr{A}_ke^{\vec a\cdot\vec
x+\vec b\cdot\vec y}\qquad\for\;k\in\mbb{N}.\eqno(2.1)$$
 According to (1.6),
$$(A_{n+1,i}-a_i)(fe^{\vec a\cdot\vec x+\vec b\cdot\vec
y})=\ptl_{x_i}(f)e^{\vec a\cdot\vec x+\vec b\cdot\vec
y},\;\;-(B_{i,n+1}+b_i)(fe^{\vec a\cdot\vec x+\vec b\cdot\vec
y})=\ptl_{y_i}(f)e^{\vec a\cdot\vec x+\vec b\cdot\vec y}\eqno(2.2)$$
for $i\in\ol{1,n}$. Repeatedly applying (2.2), we obtain $e^{\vec
a\cdot\vec x+\vec b\cdot\vec y}\in \msr{M}$. Equivalently,
$\msr{A}_{\vec a,\vec b,0}\subset\msr{M}$.

Suppose $\msr{A}_{\vec a,\vec b,\ell}\subset\msr{M}$ for some
$\ell\in\mbb{N}$.  Take any $ge^{\vec a\cdot\vec x+\vec b\cdot\vec
y}\in \msr{A}_{\vec a,\vec b,\ell}$. Since
$$(x_i\ptl_{x_1}-y_1\ptl_{y_i})(g)e^{\vec a\cdot\vec x+\vec b\cdot\vec
y},(y_i\ptl_{x_1}-y_1\ptl_{x_i})(g)e^{\vec a\cdot\vec x+\vec
b\cdot\vec y}\in \msr{A}_{\vec a,\vec
b,\ell}\subset\msr{M},\eqno(2.3)$$ we have
$$A_{i,1}(ge^{\vec a\cdot\vec x+\vec b\cdot\vec y})\equiv (a_1x_i-b_iy_1)ge^{\vec a\cdot\vec x+\vec b\cdot\vec y}
\equiv 0\;\;(\mbox{mod}\;\msr M)\eqno(2.4)$$ and
$$C_{i,1}(ge^{\vec a\cdot\vec x+\vec b\cdot\vec y})\equiv
(a_1y_i-a_iy_1)ge^{\vec a\cdot\vec x+\vec b\cdot\vec y} \equiv
0\;\;(\mbox{mod}\;\msr M)\eqno(2.5)$$ for $i\in\ol{1,n}$ by (1.5).
On the other hand, (1.4) implies
$$(D+c)(g)e^{\vec a\cdot\vec x+\vec
b\cdot\vec y}\in \msr{A}_{\vec a,\vec
b,\ell}\subset\msr{M},\eqno(2.6)$$ and so (1.6) gives
$$-A_{n+1,n+1}(ge^{\vec a\cdot\vec x+\vec b\cdot\vec y})\equiv
[\sum_{i=1}^n(a_ix_i+b_iy_i)]ge^{\vec a\cdot\vec x+\vec b\cdot\vec
y} \equiv 0\;\;(\mbox{mod}\;\msr M)\eqno(2.7)$$

Substituting (2.4) and (2.5) into (2.7), we
get$$(\sum_{i=1}^na_ib_i)y_1ge^{\vec a\cdot\vec x+\vec b\cdot\vec y}
\equiv 0\;\;(\mbox{mod}\;\msr M).\eqno(2.8)$$ Equivalently,
$y_1ge^{\vec a\cdot\vec x+\vec b\cdot\vec y}\in\msr M$. Substituting
it to (2.4) and (2.5), we obtain
$$x_ige^{\vec a\cdot\vec x+\vec
b\cdot\vec y},y_ige^{\vec a\cdot\vec x+\vec b\cdot\vec y}\in\msr
M\eqno(2.9)$$ for $i\in\ol{1,n}$. Therefore, $\msr{A}_{\vec a,\vec
b,\ell+1}\subset\msr{M}$. By induction, $\msr{A}_{\vec a,\vec b,
\ell}\subset\msr{M}$ for any $\ell\in\mbb{N}$. So $\msr{A}_{\vec
a,\vec b}=\msr{M}$. Hence $\msr{A}_{\vec a,\vec b}$ is an
irreducible $o(2n+2,\mbb{C})$-module. $\qquad\Box$\psp

 Fix $n_1,n_2\in\ol{1,n}$ with $n_1\leq n_2$. To make notations more distinguishable,  we write
$$D_{n_1,n_2}=-\sum_{i=1}^{n_1}x_i\ptl_{x_i}
+\sum_{r=n_1+1}^nx_r\ptl_{x_r}+\sum_{j=1}^{n_2}y_j\ptl_{y_j}-\sum_{s=n_2+1}^ny_s\ptl_{y_s},\eqno(2.10)$$
$$\eta_{n_1,n_2}=\sum_{i=1}^{n_1}y_i\ptl_{x_i}+\sum_{r=n_1+1}^{n_2}x_ry_r+\sum_{s=n_2+1}^n
x_s\ptl_{y_s}\eqno(2.11)$$and
$$\td c=c+n_2-n_1-n.\eqno(2.12)$$
Then we have the following representation $\pi_c^{n_1,n_2}$ of the
Lie algebra $o(2n+2,\mbb{C})$  determined by
$$\pi_c^{n_1,n_2}(A_{i,j})=E_{i,j}^x-E_{j,i}^y\eqno(2.13)$$ with
$$E_{i,j}^x=\left\{\begin{array}{ll}-x_j\ptl_{x_i}-\delta_{i,j}&\mbox{if}\;
i,j\in\ol{1,n_1},\\ \ptl_{x_i}\ptl_{x_j}&\mbox{if}\;i\in\ol{1,n_1},\;j\in\ol{n_1+1,n},\\
-x_ix_j &\mbox{if}\;i\in\ol{n_1+1,n},\;j\in\ol{1,n_1},\\
x_i\partial_{x_j}&\mbox{if}\;i,j\in\ol{n_1+1,n}
\end{array}\right.\eqno(2.14)$$
and
$$E_{i,j}^y=\left\{\begin{array}{ll}y_i\ptl_{y_j}&\mbox{if}\;
i,j\in\ol{1,n_2},\\ -y_iy_j&\mbox{if}\;i\in\ol{1,n_2},\;j\in\ol{n_2+1,n},\\
\ptl_{y_i}\ptl_{y_j} &\mbox{if}\;i\in\ol{n_2+1,n},\;j\in\ol{1,n_2},\\
-y_j\partial_{y_i}-\delta_{i,j}&\mbox{if}\;i,j\in\ol{n_2+1,n},
\end{array}\right.\eqno(2.15)$$
and
$$\pi_c^{n_1,n_2}(E_{i,n+1+j})=\left\{\begin{array}{ll}
\ptl_{x_i}\ptl_{y_j}&\mbox{if}\;i\in\ol{1,n_1},\;j\in\ol{1,n_2},\\
-y_j\ptl_{x_i}&\mbox{if}\;i\in\ol{1,n_1},\;j\in\ol{n_2+1,n},\\
x_i\ptl_{y_j}&\mbox{if}\;i\in\ol{n_1+1,n},\;j\in\ol{1,n_2},\\
-x_iy_j&\mbox{if}\;i\in\ol{n_1+1,n},\;j\in\ol{n_2+1,n},\end{array}\right.\eqno(2.16)$$
$$\pi_c^{n_1,n_2}(E_{n+1+i,j})=\left\{\begin{array}{ll}
-x_jy_i&\mbox{if}\;j\in\ol{1,n_1},\;i\in\ol{1,n_2},\\
-x_j\ptl_{y_i}&\mbox{if}\;j\in\ol{1,n_1},\;i\in\ol{n_2+1,n},\\
y_i\ptl_{x_j}&\mbox{if}\;j\in\ol{n_1+1,n},\;i\in\ol{1,n_2},\\
\ptl_{x_j}\ptl_{y_i}&\mbox{if}\;j\in\ol{n_1+1,n},\;i\in\ol{n_2+1,n},\end{array}\right.\eqno(2.17)$$
$$\pi_c^{n_1,n_2}(A_{n+1,n+1})=-D_{n_1,n_2}-\td c,\eqno(2.18)$$
$$\pi_c^{n_1,n_2}(A_{n+1,i})=\left\{\begin{array}{ll}-x_i&\mbox{if}\;\;i\in\ol{1, n_1},\\
\ptl_{x_i}&\mbox{if}\;\;i\in\ol{n_1+1,n},\end{array}\right.\eqno(2.19)$$
$$\pi_c^{n_1,n_2}(B_{i,n+1})=\left\{\begin{array}{ll}-\ptl_{y_i}&\mbox{if}\;\;\in\ol{1, n_2},\\
y_i&\mbox{if}\;\;i\in\ol{n_2+1,n},\end{array}\right. \eqno(2.20)$$
$$\pi_c^{n_1,n_2}(A_{i,n+1})=\left\{\begin{array}{ll} \eta_{n_1,n_2}\ptl_{y_i}-(D_{n_1,n_2}+\td c-1)\ptl_{x_i}&\mbox{if}\;\;i\in\ol{1, n_1},\\
 \eta_{n_1,n_2}\ptl_{y_i}-x_i(D_{n_1,n_2}+\td c)&\mbox{if}\;\;i\in\ol{n_1+1, n_2},\\
 - \eta_{n_1,n_2}y_i-x_i(D_{n_1,n_2}+\td
 c)&\mbox{if}\;\;i\in\ol{n_2+1,n},\end{array}\right.\eqno(2.21)$$
$$\pi_c^{n_1,n_2}(C_{n+1,i})
=\left\{\begin{array}{ll} \eta_{n_1,n_2} x_i+y_i(D_{n_1,n_2}+\td c)&\mbox{if}\;\;i\in\ol{1, n_1},\\
- \eta_{n_1,n_2}\ptl_{x_i}+y_i(D_{n_1,n_2}+\td c)&\mbox{if}\;\;i\in\ol{n_1+1, n_2},\\
 -\eta_{n_1,n_2}\ptl_{x_i}+(D_{n_1,n_2}+\td
 c-1)\ptl_{y_i}&\mbox{if}\;\;i\in\ol{n_2+1,n}\end{array}\right.\eqno(2.22)$$
for $i,j\in\ol{1,n}$.

Set
$$\msr A_{\la k\ra}=\mbox{Span}\{x^\al
y^\be\mid\al,\be\in\mbb{N}\:^n;\sum_{r=n_1+1}^n\al_r-\sum_{i=1}^{n_1}\al_i+
\sum_{i=1}^{n_2}\be_i-\sum_{r=n_2+1}^n\be_r=k\}\eqno(2.23)$$ for
$k\in\mbb{Z}$. Then
$$\msr A_{\la k\ra}=\{u\in\msr A\mid D_{n_1,n_2}(u)=k
u\}\eqno(2.24)$$ Observe that the Lie subalgebra
$$\msr
K=\sum_{i,j=1}^n(\mbb{C}A_{i,j}+\mbb{C}B_{i,j}+\mbb{C}C_{i,j})\cong
o(2n,\mbb{C}).\eqno(2.25)$$ With respect to the presentation
$\pi_c^{n_1,n_2}$, $\msr A_{\la k\ra}$ forms a $\msr K$-module.
Write
$$\msr D_{n_1,n_2}=-\sum_{i=1}^{n_1}x_i\ptl_{y_i}+\sum_{r=n_1+1}^{n_2}\ptl_{x_r}\ptl_{y_r}-\sum_{s=n_2+1}^n
y_s\ptl_{x_s}.\eqno(2.26)$$ Note that
 as operators on $\msr A$,
$$\xi\eta_{n_1,n_2}=\eta_{n_1,n_2}\xi,\;\;\xi\msr D_{n_1,n_2} =
\msr D_{n_1,n_2}\xi\qquad\for\;\;\xi\in\msr K.\eqno(2.27)$$ In
particular,
$$\msr H_{\la k\ra}=\{u\in\msr A_{\la k\ra}\mid \msr D_{n_1,n_2}(u)=0
\}\eqno(2.28)$$ forms a $\msr K$-module. The following result is
taken from Luo and the author's work [LX2].\pse

{\bf Lemma 2.2}. {\it For any $n_1-n_2+1-\dlt_{n_1,n_2}\geq
k\in\mbb{Z}$, $\msr H_{\la k\ra}$ is an irreducible $\msr
K$-submodule and $\msr A_{\la
k\ra}=\bigoplus_{i=0}^\infty\eta_{n_1,n_2}^i(\msr H_{\la k-2i\ra})$
is a decomposition of irreducible $\msr K$-submodules.}\psp

Now we have the second result in this section.\psp

{\bf Theorem 2.3}. {\it The representation $\pi_c^{n_1,n_2}$ of
$o(2n+2,\mbb{C})$ on $\msr A$ is irreducible if $c\not\in
\mbb{Z}/2$.}

{\it Proof}. Let $\msr M$ be a nonzero $o(2n+2,\mbb{C})$-submodule
of $\msr A$. By (2.18) and (2.24),
$$\msr M=\bigoplus_{k\in\mbb{Z}}\msr A_{\la k\ra}\bigcap
\msr M.\eqno(2.29)$$ Thus $\msr A_{\la k\ra}\bigcap \msr M\neq\{0\}$
for some $k\in \mbb{Z}$. If $k>n_1-n_2+1-\dlt_{n_1,n_2}$, then
$$
\{0\}\neq(-x_1)^{k-(n_1-n_2+1-\dlt_{n_1,n_2})}(\msr A_{\la
k\ra}\bigcap \msr M) =A_{n+1,1}^{k-(n_1-n_2+1-\dlt_{n_1,n_2})}(\msr
A_{\la k\ra}\bigcap \msr M)\eqno(2.30)$$ by (2.19),  which implies
$\msr A _{\la n_1-n_2+1-\dlt_{n_1,n_2} \ra}\bigcap \msr M\neq
\{0\}$. Thus we can assume $k\leq n_1-n_2+1-\dlt_{n_1,n_2}$. Observe
that the Lie subalgebra
$$\msr L=\sum_{i,j=1}^n\mbb CA_{i,j}\cong
sl(n,\mbb{C}).\eqno(2.31)$$ By Lemma 2.2, $\msr A_{\la
k\ra}=\bigoplus_{i=0}^\infty\eta_{n_1,n_2}^i(\msr H_{\la k-2i\ra})$
is a decomposition of irreducible $\msr K$-submodules. Moreover,
$\eta_{n_1,n_2}^i(\msr H_{\la k-2i\ra})$ are highest-weight $\msr
L$-modules with distinct highest weights by [LX1]. Hence
$$\eta_{n_1,n_2}^i(\msr H_{\la k-2i\ra})\subset \msr M\;\;\mbox{for some}\;\;i\in\mbb{N}.\eqno(2.32)$$

Observe that
$$x_1^{-k+2i}\in \msr H_{\la k-2i\ra}.\eqno(2.33)$$
By (2.11) and (2.20),
$$i!(-1)^i(\prod_{r=1}^i(-k+i+r))x_1^{-k+i}=B_{1,2n+2}^i(\eta^i_{n_1,n_2}(x_1^{-k+2i}))\in
\msr M.\eqno(2.34)$$ Thus
$$\msr H_{\la k-i\ra}\subset \msr M.\eqno(2.35)$$
So we can just assume
$$\msr H_{\la k\ra}\subset \msr M.\eqno(2.36)$$
According to (2.19),
$$x_1^{-k+s}=(-1)^sA_{n+1,1}^s(x_1^{-k})\in
\msr M\qquad\for\;\;s\in\mbb{N}.\eqno(2.37)$$ So Lemma 2.2 gives
$$\msr H_{\la k-s\ra}\subset \msr M\qquad\for\;\;s\in\mbb{N}.\eqno(2.38)$$
For any $r\in k-\mbb{N}$, we suppose
$\eta_{n_1,n_2}^s(x_1^{-r+s}),\eta_{n_1,n_2}^s(x_1^{-r+s+1})\in \msr
M$ for some $s\in\mbb N$. Applying (2.22) to it, we get
$$C_{n+1,1}[\eta_{n_1,n_2}^s(x_1^{-r+s})]=\eta_{n_1,n_2}^{s+1}(x_1^{-r+s+1})
+(r+\td c)\eta_{n_1,n_2}^s(y_1x_1^{-r+s}) \in\msr M.\eqno(2.39)$$ By
(2.11) and (2.22),
$$C_{n+1,i}[\eta_{n_1,n_2}^s(x_1^{-r+s+1})]=(r-1+\td
c)\eta_{n_1,n_2}^s(y_ix_1^{-r+s+1})\in\msr M\eqno(2.40)$$ for
$i\in\ol{n_1+1,n_2}$. According to (2.11) and (2.21),
$$A_{i,n+1}[\eta_{n_1,n_2}^s(y_ix_1^{-r+s+1})]=\eta_{n_1,n_2}^{s+1}(x_1^{-r+s+1})
-(r+\td c)\eta_{n_1,n_2}^s(x_iy_ix_1^{-r+s+1})\in \msr
M\eqno(2.41)$$ for $i\in\ol{n_1+1,n_2}$. Again (2.11), (2.39) and
(2.41) lead to
$$ (1+r+\td
c-n_2+n_1)\eta_{n_1,n_2}^{s+1}(x_1^{-r+s+1})\in\msr{M}\Rightarrow
\eta_{n_1,n_2}^{s+1}(x_1^{-r+s+1})\in\msr{M}.\eqno(2.42)$$

By induction,
$$\eta_{n_1,n_2}^\ell(x_1^{-r+\ell})\in\msr{M}\qquad\for\;\;\ell\in\mbb
N.\eqno(2.43)$$ Since $\eta_{n_1,n_2}^\ell(\msr H_{\la
r-\ell\ra})\ni \eta_{n_1,n_2}^\ell(x_1^{-r+\ell})$ is an irreducible
$\msr L$-module by Lemma 2.2, we have
$$\eta_{n_1,n_2}^\ell(\msr H_{\la
r-\ell\ra})\subset\msr M \qquad\for\;\;\ell\in\mbb N.\eqno(2.44)$$
Taking $r=m-\ell$ with $m\in k-\mbb N$, we get
$$\eta_{n_1,n_2}^\ell(\msr H_{\la
m-2\ell\ra})\subset\msr M \qquad\for\;\;\ell\in\mbb N.\eqno(2.45)$$
According to Lemma 2.2,
$$\msr A_{\la
m\ra}=\bigoplus_{\ell=0}^\infty\eta_{n_1,n_2}^\ell(\msr H_{\la
m-2\ell\ra})\subset\msr M\qquad\for\;\;m\in k-\mbb N.\eqno(2.46)$$

Expression (2.21) gives
$$\pi_c^{n_1,n_2}(A_{i,n+1})y_i=\left\{\begin{array}{ll} \eta_{n_1,n_2}(y_i\ptl_{y_i}+1)-y_i\ptl_{x_i}(D_{n_1,n_2}+\td c+1)&\mbox{if}\;\;i\in\ol{1, n_1},\\
 \eta_{n_1,n_2}(y_i\ptl_{y_i}+1)-x_iy_i(D_{n_1,n_2}+\td c+1)&\mbox{if}\;\;i\in\ol{n_1+1,n_2},\end{array}\right.
 \eqno(2.47)$$
$$\pi_c^{n_1,n_2}(A_{j,n+1})\ptl_{y_j}=- \eta_{n_1,n_2}y_j\ptl_{y_j}-x_j\ptl_{y_j}(D_{n_1,n_2}+\td
 c+1)\qquad\for\;\;j\in\ol{n_2+1,n}.\eqno(2.48)$$
Moreover, (2.22) yields
$$\pi_c^{n_1,n_2}(C_{n+1,r})\ptl_{x_r}=\eta_{n_1,n_2}x_r\ptl_{x_i} +y_r\ptl_{x_r}(D_{n_1,n_2}+\td
c+1)\qquad\for\;\;r\in\ol{1,n_1},\eqno(2.49)$$
\begin{eqnarray*}& &\pi_c^{n_1,n_2}(C_{n+1,s})x_s \\&=&\left\{\begin{array}{ll}
- \eta_{n_1,n_2}(x_s\ptl_{x_s}+1)+x_sy_s(D_{n_1,n_2}+\td c+1)&\mbox{if}\;\;s\in\ol{n_1+1,n_2},\\
 -\eta_{n_1,n_2}(x_s\ptl_{x_s}+1)+x_s\ptl_{y_s}(D_{n_1,n_2}+\td
 c+1)&\mbox{if}\;\;s\in\ol{n_2+1,n}.\end{array}\right.\hspace{2.3cm}(2.50)\end{eqnarray*}
Thus
\begin{eqnarray*}\hspace{2cm}& &\sum_{i=1}^{n_2}\pi_c^{n_1,n_2}(A_{i,n+1})y_i+\sum_{j=n_2+1}^n\pi_c^{n_1,n_2}(A_{j,n+1})\ptl_{y_j}
\\
&&-\sum_{r=1}^{n_1}\pi_c^{n_1,n_2}(C_{n+1,r})\ptl_{x_r}-\sum_{s=n_1+1}^n\pi_c^{n_1,n_2}(C_{n+1,s})x_s\\
&=&\eta_{n_1,n_2}(-D_{n_1,n_2}+n_2+n-n_1-2(\td
c+1))\hspace{4.7cm}(2.51)\end{eqnarray*} as operators on $\msr A$.
Suppose that $\msr A_{\la \ell-s\ra}\subset \msr M$ for some
$k\leq\ell\in\mbb{Z}$ and any $s\in\mbb{N}$. For any $f\in \msr
A_{\la \ell-1\ra}$, we apply the above equation to it and get
$$(1-\ell+n_2+n-n_1-2(\td
c+1))\eta_{n_1,n_2}(f)\in \msr M.\eqno(2.52)$$ Since $c\not\in \mbb
Z/2$, we have
$$\eta_{n_1,n_2}(f)\in \msr M.\eqno(2.53)$$

Now for any $g\in\msr A_{\la \ell\ra}$, we have $\ptl_{y_1}(g)\in
\msr A_{\la \ell-1\ra}$. By (2.21),
$$A_{1,n+1}(g)=\eta_{n_1,n_2}(\ptl_{y_1}(g))-(\ell+\td
c)\ptl_{x_1}(g)\in\msr M.\eqno(2.54)$$ Moreover, (2.53) and (2.54)
yield
$$\ptl_{x_1}(g)\in \msr M\qquad\for\;\;g\in \msr A_{\la
\ell\ra}.\eqno(2.55)$$ Since
$$\ptl_{x_1}(\msr A_{\la
\ell\ra})=\msr A_{\la \ell+1\ra},\eqno(2.56)$$ we obtain
$$\msr A_{\la \ell+1\ra}\subset \msr M.\eqno(2.57)$$
By induction on $\ell$, we find
$$\msr A_{\la \ell\ra}\subset \msr M\qquad\for\;\;\ell\in\mbb{Z},\eqno(2.58)$$
or equivalently, $\msr A=\bigoplus_{\ell\in\mbb{Z}}\msr A_{\la
\ell\ra}=\msr M$. Thus $\msr A$ is an irreducible $o(2n+2,\mbb
C)$-module.$\qquad\Box$\psp

{\bf Remark 2.4}. The above irreducible representation depends on
the three parameters $c\in \mbb{F}$ and $m_1,m_2\in\ol{1,n}$. It is
not highest-weight type because of the mixture of multiplication
operators and differential operators in
 (2.16), (2.17) and
(2.19)-(2.22). Since $\msr A$ is not completely reducible as a $\msr
L$-module by [LX1] when $n\geq 2$ and $n_1<n$, $\msr A$  is not a
unitary $o(2n+2,\mbb{C})$-module. Expression (2.18) shows that $\msr
A$ is a weight $o(2n+2,\mbb{C})$-module with finite-dimensional
weight subspaces.\psp

Theorem 1 follows from Theorem 2.1, Theorem 2.3 and the above
remark.

\section{Proof of Theorem 2}

$\quad\;$ In this section, we prove Theorem 2. Our first result in
this section is as follows.\psp

{\bf Theorem 3.1}. {\it The representation $\pi_{c,\vec a,\vec b}$
of $o(2n+3,\mbb{C})$ is irreducible for any $c\in\mbb{C}$ if
$a_0^2+2\sum_{i=1}^na_ib_i\neq 0$.}

{\it Proof}. Let $\msr B_k$ be the subspace of homogeneous
polynomials with degree $k$. Set
$$\msr B_{\vec a,\vec b,k}=\msr B_ke^{\vec a\cdot\vec
x+\vec b\cdot\vec y}\qquad\for\;k\in\mbb{N}\eqno(3.1)$$ (cf. (1.15)
and the second equation in (1.8)). Let ${\msr M}$ be a nonzero
$o(2n+3,\mbb{C})$-submodule of $\msr B_{\vec a,\vec b}$. Take any
$0\neq fe^{\vec a\cdot\vec x+\vec b\cdot\vec y}\in \msr{M}$ with
$f\in \msr B$.  According to (1.6),
$$(A_{n+1,i}-a_i)(fe^{\vec a\cdot\vec x+\vec b\cdot\vec
y})=\ptl_{x_i}(f)e^{\vec a\cdot\vec x+\vec b\cdot\vec
y},\;\;-(B_{i,n+1}+b_i)(fe^{\vec a\cdot\vec x+\vec b\cdot\vec
y})=\ptl_{y_i}(f)e^{\vec a\cdot\vec x+\vec b\cdot\vec y}\eqno(3.2)$$
for $i\in\ol{1,n}$. Moreover, the second equation in (1.14) gives
$$-(K_{2n+2}+a_0)(fe^{\vec a\cdot\vec x+\vec b\cdot\vec
y})=\ptl_{x_0}(f)e^{\vec a\cdot\vec x+\vec b\cdot\vec
y}.\eqno(3.3)$$
 Repeatedly applying (3.2) and (3.3), we obtain $e^{\vec
a\cdot\vec x+\vec b\cdot\vec y}\in \msr{M}$. Equivalently, $\msr
B_{\vec a,\vec b,0}\subset\msr{M}$. Suppose $\msr B_{\vec a,\vec
b,\ell}\subset\msr{M}$ for some $\ell\in\mbb{N}$. Let $ ge^{\vec
a\cdot\vec x+\vec b\cdot\vec y}$ be any element in $\msr{A}_{\vec
a,\vec b,\ell}$.\pse

{\it Case 1. $a_i\neq 0$ or $b_i\neq 0$ for some
$i\in\ol{1,n}$.}\pse

 By symmetry, we may assume $a_1\neq 0$. Expression (2.3) with $\msr A_{\vec a,\vec
 b,\ell}$ replaced by $\msr B_{\vec a,\vec b,\ell}$ implies
$$A_{i,1}(ge^{\vec a\cdot\vec x+\vec b\cdot\vec y})\equiv (a_1x_i-b_iy_1)ge^{\vec a\cdot\vec x+\vec b\cdot\vec y}
\equiv 0\;\;(\mbox{mod}\;\msr M)\eqno(3.4)$$ and
$$C_{1+i,1}(ge^{\vec a\cdot\vec x+\vec b\cdot\vec y})\equiv
(a_1y_i-a_iy_1)ge^{\vec a\cdot\vec x+\vec b\cdot\vec y} \equiv
0\;\;(\mbox{mod}\;\msr M)\eqno(3.5)$$ for $i\in\ol{1,n}$ by (1.5).
Moreover, the first equation in (1.13) gives
$$K_1(ge^{\vec a\cdot\vec x+\vec b\cdot\vec y})\equiv
(a_1x_0-a_0y_1)ge^{\vec a\cdot\vec x+\vec b\cdot\vec y} \equiv
0\;\;(\mbox{mod}\;\msr M)\eqno(3.6)$$ because
$$(x_0\ptl_{x_1}-y_1\ptl_{x_0})(g)e^{\vec a\cdot\vec x+\vec b\cdot\vec
y}\in\msr B_{\vec a,\vec b,\ell}\subset\msr{M}.\eqno(3.7)$$

 On the other hand, the
second equation in (1.6) with $D$ in (1.12) gives
$$-A_{n+1,n+1}(ge^{\vec a\cdot\vec x+\vec b\cdot\vec y})\equiv
[a_0x_0+\sum_{i=1}^n(a_ix_i+b_iy_i)]ge^{\vec a\cdot\vec x+\vec
b\cdot\vec y} \equiv 0\;\;(\mbox{mod}\;\msr M)\eqno(3.8)$$ by (2.6)
with $\msr A_{\vec a,\vec
 b,\ell}$ replaced by $\msr B_{\vec a,\vec b,\ell}$.
 Substituting (3.4)-(3.6) into (3.8), we get
$$(a_0^2+2\sum_{i=1}^na_ib_i)y_1ge^{\vec
a\cdot\vec x+\vec b\cdot\vec y} \equiv 0\;\;(\mbox{mod}\;\msr
M).\eqno(3.9)$$ Equivalently, $y_1ge^{\vec a\cdot\vec x+\vec
b\cdot\vec y}\in\msr M$. Substituting it to (3.4)-(3.6), we obtain
$$x_0ge^{\vec a\cdot\vec x+\vec
b\cdot\vec y},x_ige^{\vec a\cdot\vec x+\vec b\cdot\vec
y},y_ige^{\vec a\cdot\vec x+\vec b\cdot\vec y}\in\msr M\eqno(3.10)$$
for $i\in\ol{1,n}$. Therefore, $\msr B_{\vec a,\vec
b,\ell+1}\subset\msr{M}$. By induction, $\msr B_{\vec a,\vec b,
\ell}\subset\msr{M}$ for any $\ell\in\mbb{N}$. So $\msr B_{\vec
a,\vec b}=\msr{M}$. Hence $\msr B_{\vec a,\vec b}$ is an irreducible
$o(2n+3,\mbb{C})$-module.\pse

{\it Case 2. $a_0\neq 0$ and $a_i=b_0=0$ for $i\in\ol{1,n}$.}\pse

Under the above assumption,
$$K_i(ge^{\vec a\cdot\vec x+\vec b\cdot\vec
y})=(x_0\ptl_{x_i}-y_i\ptl_{x_0}-a_0y_i)(g)e^{\vec a\cdot\vec x+\vec
b\cdot\vec y}\in\msr M\eqno(3.11)$$ and
$$K_{n+1+i}(ge^{\vec a\cdot\vec x+\vec b\cdot\vec
y})=(x_0\ptl_{y_i}-x_i\ptl_{x_0}-a_0x_i)(g)e^{\vec a\cdot\vec x+\vec
b\cdot\vec y}\in\msr M\eqno(3.12)$$ for $i\in\ol{1,n}$. Note
$$(x_0\ptl_{x_i}-y_i\ptl_{x_0})(g)e^{\vec a\cdot\vec x+\vec
b\cdot\vec y}, (x_0\ptl_{y_i}-x_i\ptl_{x_0})(g)e^{\vec a\cdot\vec
x+\vec b\cdot\vec y}\in\msr B_{\vec a,\vec b,\ell}\subset\msr{M}
\eqno(3.13)$$ by the inductional assumption. Thus (3.10) and (3.11)
imply
$$y_ige^{\vec a\cdot\vec x+\vec
b\cdot\vec y},x_ige^{\vec a\cdot\vec x+\vec b\cdot\vec y}\in\msr
M\qquad\for\;\;i\in\ol{1,n}.\eqno(3.14)$$ Now (3.8) yields
$x_0ge^{\vec a\cdot\vec x+\vec b\cdot\vec y}\in\msr M$. So  $B_{\vec
a,\vec b,\ell+1}\subset\msr{M}$. By induction, $\msr B=\msr M$; that
is, $\msr B$ is irreducible. $\qquad\Box$\psp

Fix $n_1,n_2\in\ol{1,n}$ with $n_1\leq n_2$ . Reset
$$D_{n_1,n_2}=x_0\ptl_{x_0}-\sum_{i=1}^{n_1}x_i\ptl_{x_i}
+\sum_{r=n_1+1}^nx_r\ptl_{x_r}+\sum_{j=1}^{n_2}y_j\ptl_{y_j}-\sum_{s=n_2+1}^ny_s\ptl_{y_s},\eqno(3.15)$$
$$\msr D_{n_1,n_2}=\ptl_{x_0}^2-2\sum_{i=1}^{n_1}x_i\ptl_{y_i}+2\sum_{r=n_1+1}^{n_2}\ptl_{x_r}\ptl_{y_r}-2\sum_{s=n_2+1}^n
y_s\ptl_{x_s}\eqno(3.16)$$ and
$$\eta_{n_1,n_2}=\frac{x_0^2}{2}+\sum_{i=1}^{n_1}y_i\ptl_{x_i}+\sum_{r=n_1+1}^{n_2}x_ry_r+\sum_{s=n_2+1}^n
x_s\ptl_{y_s}.\eqno(3.17)$$ Then the representation
$\pi_c^{n_1,n_2}$ of $o(2n+3,\mbb C)$ is determined as follows:
$\pi_c|_{o(2n+2,\mbb C)}$ is given by (2.12)-(2.22) with
$D_{n_1,n_2}$ in (3.15) and $\eta_{n_1,n_2}$ in (3.17), and
$$\pi_c^{n_1,n_2}(K_i)=\left\{\begin{array}{ll}-x_0x_i-y_i\ptl_{x_0}&\mbox{if}\;i\in\ol{1,n_1},\\
x_0\ptl_{x_i}-y_i\ptl_{x_0}&\mbox{if}\;i\in\ol{n_1+1,n_2},\\
x_0\ptl_{x_i}-\ptl_{x_0}\ptl_{y_i}&\mbox{if}\;i\in\ol{n_2+1,n},\end{array}\right.\eqno(3.18)$$
$$\pi_c^{n_1,n_2}(K_{n+1+i})=\left\{\begin{array}{ll}x_0\ptl_{y_i}-\ptl_{x_0}\ptl_{x_i}&\mbox{if}\;i\in\ol{1,n_1},\\
x_0\ptl_{y_i}-x_i\ptl_{x_0}&\mbox{if}\;i\in\ol{n_1+1,n_2},\\
-x_0y_i-x_i\ptl_{x_0}&\mbox{if}\;i\in\ol{n_2+1,n},\end{array}\right.\eqno(3.19)$$
$$\pi_c^{n_1,n_2}(K_{n+1})=x_0(D_{n_1,n_2}+\td c)-\eta_{n_1,n_2}\ptl_{x_0},\qquad \pi_c^{n_1,n_2}(K_{2n+2})=-\ptl_{x_0}.
\eqno(3.20)$$ Note that
$$\msr G=\msr K+\sum_{i=1}^{2n+2}\mbb{C}K_i\eqno(3.21)$$
is a Lie subalgebra isomorphic to $o(2n+1,\mbb C)$.

Define
$$\msr B_{\la k\ra}=\sum_{i=}^\infty\msr A_{\la
k\ra}x_0^i.\eqno(3.22)$$ Then
$$\msr B_{\la k\ra}=\{u\in\msr B\mid D_{n_1,n_2}(u)=k
u\}\qquad\for\;\;k\in\mbb Z\eqno(3.23)$$ and
$$\msr B=\bigoplus_{k\in\mbb Z}\msr B_{\la k\ra}.\eqno(3.24)$$
Moreover,
$$\xi D_{n_1,n_2} =
D_{n_1,n_2}\xi,\;\;\xi\eta_{n_1,n_2}=\eta_{n_1,n_2}\xi,\;\;\xi\msr
D_{n_1,n_2} = \msr D_{n_1,n_2}\xi\qquad\for\;\;\xi\in\msr
G\eqno(3.25)$$ as operators on $\msr B$. In particular, $\msr B_{\la
k\ra}$ forms a $\msr G$-module for any $k\in\msr Z$. Furthermore,
$$\msr H_{\la k\ra}=\{u\in\msr B_{\la k\ra}\mid \msr D_{n_1,n_2}(u)=0
\}\eqno(3.26)$$ forms a $\msr G$-module. The following result is
taken from Luo and the author's work [LX2].\pse

{\bf Lemma 3.2}. {\it For any $ k\in\mbb{Z}$, $\msr H_{\la k\ra}$ is
an irreducible $\msr G$-submodule and $\msr A_{\la
k\ra}=\bigoplus_{i=0}^\infty\eta_{n_1,n_2}^i(\msr H_{\la k-2i\ra})$
is a decomposition of irreducible $\msr G$-submodules.}\psp

Now we have the second result in this section.\psp

{\bf Theorem 3.3}. {\it The representation $\pi_c^{n_1,n_2}$ of
$o(2n+3,\mbb{C})$ on $\msr A$ is irreducible if $c\not\in
\mbb{Z}/2$.}

{\it Proof}. Let $\msr M$ be a nonzero $o(2n+3,\mbb{C})$-submodule
of $\msr B$. By (3.23) and (2.18) with $D_{n_1,n_2}$ in (3.15),
$$\msr M=\bigoplus_{k\in\mbb{Z}}\msr B_{\la k\ra}\bigcap
\msr M.\eqno(3.27)$$ Thus $\msr B_{\la k\ra}\bigcap \msr M\neq\{0\}$
for some $k\in \mbb{Z}$.  Take the Lie subalgebra $\msr L$ in
(2.31). By Lemma 3.2, $\msr B_{\la
k\ra}=\bigoplus_{i=0}^\infty\eta_{n_1,n_2}^i(\msr H_{\la k-2i\ra})$
is a decomposition of irreducible $\msr G$-submodules. Moreover,
$\eta_{n_1,n_2}^i(\msr H_{\la k-2i\ra})$ are highest-weight $\msr
L$-modules with distinct highest weights by [LX1]. Hence
$$\eta_{n_1,n_2}^i(\msr H_{\la k-2i\ra})\subset \msr M\;\;\mbox{for some}\;\;i\in\mbb{N}.\eqno(3.28)$$
Lemma 3.2 and the arguments in (2.33)-(2.36) show
$$\msr H_{\la k-s\ra}\subset \msr M\qquad\for\;\;s\in\mbb{N}.\eqno(3.29)$$

Suppose $\eta_{n_1,n_2}^s(x_1^{-r+s})\in \msr M$ for any $r\in
k-\mbb N$ and some $s\in\mbb N$. Then,
$$K_{n+1}(\eta_{n_1,n_2}^s(x_1^{-r+s+1}))=(r-1+\td c)\eta_{n_1,n_2}^s(x_0x_1^{-r+s+1})\in\msr
M\eqno(3.30)$$ by (3.17) and the first equation in (3.20), which
implies $\eta_{n_1,n_2}^s(x_0x_1^{-r+s+1})\in\msr M.$ Moreover,
$$K_{n+1}(\eta_{n_1,n_2}^s(x_0x_1^{-r+s+1}))=-\eta_{n_1,n_2}^{s+1}(x_1^{-r+s+1})+(r+\td c)
\eta_{n_1,n_2}^s(x_0^2x_1^{-r+s+1}).\eqno(3.31)$$ Now (2.39) and
(2.41) with $\eta_{n_1,n_2}$ in (3.17),  and (3.31) lead to
$$ (1/2+r+\td
c-n_2+n_1)\eta_{n_1,n_2}^{s+1}(x_1^{-r+s+1})\in\msr{M}\Rightarrow
\eta_{n_1,n_2}^{s+1}(x_1^{-r+s+1})\in\msr{M}.\eqno(3.32)$$ By
induction,
$$\eta_{n_1,n_2}^\ell(x_1^{-r+\ell})\in\msr{M}\qquad\for\;\;\ell\in\mbb
N,\;r\in k-\mbb N.\eqno(3.33)$$  According to Lemma 3.2,
$$\msr B_{\la
m\ra}=\bigoplus_{\ell=0}^\infty\eta_{n_1,n_2}^\ell(\msr H_{\la
m-2\ell\ra})\subset\msr M\qquad\for\;\;m\in k-\mbb N.\eqno(3.34)$$

Observe that
$$\pi_c^{n_1,n_2}(K_{n+1})x_0=x_0^2(D_{n_1,n_2}+\td c+1)-\eta_{n_1,n_2}(x_0\ptl_{x_0}+1)\eqno(3.35)$$
by (3.20). Then (3.35) and (2.47)-(2.50) with $\eta_{n_1,n_2}$ in
(3.17) and $D_{n_1,n_2}$ in (3.15) yield
\begin{eqnarray*}\hspace{2cm}& &-\pi_c^{n_1,n_2}(K_{n+1})x_0+\sum_{i=1}^{n_2}\pi_c^{n_1,n_2}(A_{i,n+1})y_i+\sum_{j=n_2+1}^n\pi_c^{n_1,n_2}(A_{j,n+1})\ptl_{y_j}
\\
&&-\sum_{r=1}^{n_1}\pi_c^{n_1,n_2}(C_{n+1,r})\ptl_{x_r}-\sum_{s=n_1+1}^n\pi_c^{n_1,n_2}(C_{n+1,s})x_s\\
&=&\eta_{n_1,n_2}(1-D_{n_1,n+2}+n_2+n-n_1-2(\td
c+1))\hspace{4.1cm}(3.36)\end{eqnarray*} as operators on $\msr B$.
The arguments in (2.52)-(2.58) show $\msr M=\msr B$; that is, $\msr
B$ is an irreducible $o(2n+3,\mbb C)$-module. $\qquad\Box$ \psp

{\bf Remark 3.4}. The above irreducible representation depends on
the three parameters $c\in \mbb{C}$ and $m_1,m_2\in\ol{1,n}$. It is
not highest-weight type because of the mixture of multiplication
operators and differential operators in
 (2.16), (2.17),
(2.19)-(2.22), (3.18) and (3.19). Since $\msr B$ is not completely
reducible as a $\msr L$-module by [LX1] when $n\geq 2$ and $n_1<n$,
$\msr B$ is not a unitary $o(2n+3,\mbb{C})$-module. Expression
(2.18) with $D_{n_1,n_2}$ in (3.15) shows that $\msr B$ is a weight
$o(2n+2,\mbb{C})$-module with finite-dimensional weight
subspaces.\psp

Theorem 2 follows from Theorem 3.1, Theorem 3.3 and the above
remark.

 \end{document}